\documentclass[10pt,twoside]{article}
\usepackage{graphicx}
\usepackage{amsmath}
\usepackage{Latex-document}
\usepackage{latexsym, amsfonts,amssymb}

\title{\bf  Characteristic Classes of Flat Bundles \vskip -2mm and
Determinant of the \vskip -2mm Gauss-Manin Connection \vskip 6mm}
\author{H\'{e}l\`{e}ne Esnault\vspace*{-0.5cm}\thanks{Mathematik, Universit\"{a}t Essen, FB6, Mathematik, 45117
Essen, Germany. E-mail: esnault@uni-essen.de}}
\date{\vspace{-8mm}}

\begin{document}

\maketitle

\thispagestyle{first} \setcounter{page}{471}

\begin{quotation}
\vskip 4.5mm

\end{quotation}

\vskip 12mm

\section{Introduction} \label{section 1}\setzero

\vskip-5mm \hspace{5mm}

The purpose of this note is to give a survey on recent progress on characteristic classes of flat bundles, and how
they behave in a family.

\section{Characteristic classes} \label{section 2}\setzero

\vskip-5mm \hspace{5mm}

Let $X$ be a smooth algebraic variety over a field $k$. In [13] and [15], we defined the ring
\begin{eqnarray}
  AD(X) & = & \oplus_nAD^n(X) \nonumber \\
  & = & \oplus_n \mathbb{H}^n(X, \mathcal{K}_n^M\stackrel{d \log}{\longrightarrow}\Omega_{X/k}^n
  \stackrel{d}\longrightarrow\ldots\rightarrow\Omega_{X/k}^{2n-1}) \label{2.1}
\end{eqnarray}
of algebraic differential characters. Here the Zariski sheaf
$\mathcal {K}_n^M$ is the kernel of the residue map from Milnor
$K$-theory at the generic point of $X$ to Milnor $K$-theory at
codimension 1 points. More precisely, $\mathcal {K}_n^M$ satisfies
a Gersten type resolution (see [16] and [18])
\begin{eqnarray*}
  \mathcal{K}_n^M & \stackrel{\cong}{\rightarrow} & (i_{k(X), \ast}K_n^M(k(X))\stackrel{\mathrm{Res}}{\longrightarrow}
  \oplus_{x\in X^{(1)}}i_{x, \ast}K_{n-1}^M(\kappa(x))\rightarrow \\
  & & \ldots\oplus_{x\in X^{(a)}}i_{x, \ast}K_{n-a}^M(\kappa(x))\rightarrow\ldots \rightarrow\oplus_{x\in
  X^{(n)}}i_{x, \ast}K_0^M(\kappa(x))).
\end{eqnarray*}
Here $X^{(a)}$ means the free group on points in codimension $a$,
while $i_x:x\rightarrow X$ is the embedding. The map $d\log
(\{a_1, \ldots, a_n\})=d\log a_1 \wedge\cdots\wedge d\log a_n$
from $K_n^M(k(X))$ to $\Omega_{k(X)/k}^n$ carries $\mathcal
{K}_n^M$ to
$$\Omega_{X/k}^n=\mathrm{Ker}(\Omega_{k(X)}^{n}\stackrel{\mathrm{res}}{\longrightarrow}
\oplus_{x\in X^{(1)}}\Omega_{x/k}^{n-1}).$$
This defines the map $d\log:
\mathcal{K}_n^M\rightarrow\Omega_{X/k}^n$.

By the Gersten resolution, $H^n(X, \mathcal {K}_n^M)=CH^n(X)$, the
Chow group of codimension $n$ points. Thus one has a forgetful map
\begin{equation} \label{2.2}
  AD^n(X)\stackrel{\mathrm{forget}}{\longrightarrow}CH^n(X).
\end{equation}
The restriction map to the generic point Spec$(k(X))$ fulfills
\begin{eqnarray}
  AD^1(X) & \stackrel{\cong}{\rightarrow} & H^0(X, \Omega_{X/k}^1/d\log \mathcal {O}_X^{\times})\subset AD^1(k(X))
  \nonumber \\
  AD^n(X) & {\rightarrow} & H^0(X, \Omega_{X/k}^{2n-1}/d\Omega_{X/k}^{2n-2})\subset
\Omega_{k(X)/k}^{2n-1}/d\Omega_{k(X)/k}^{2n-2}, \quad \text{for } n\geq2 \label{2.3}
\end{eqnarray}
(see [2]). It is no longer injective for $n\geq 2$.

The K\"{a}hler differential
$d:\Omega_{X/k}^{2n-1}\rightarrow\Omega_{X/k, \rm{clsd}}^{2n}$
defines
\begin{equation} \label{2.4}
  AD^n(X)\stackrel{d}{\rightarrow}H^0(X, \Omega_{X/k, \mathrm{clsd}}^{2n}).
\end{equation}

The ring $AD(X)=\oplus_nAD^n(X)$ contains the subring\\
\begin{eqnarray}
  AD(X)_{\rm{clsd}}(X) & = & \oplus_nAD_{\rm{clsd}}^n(X) \oplus_n\mathbb{H}^n(X, \mathcal {K}_n^M
  \stackrel{d \log}{\longrightarrow} \Omega_{X/k}^n\rightarrow \ldots\rightarrow\Omega_{X/k}^{\mathrm{dim}(X)})
  \nonumber \\
  & = & \mathrm{Ker}(AD(X)\stackrel{d}{\rightarrow}\oplus_nH^0(X, \Omega_{X/k, \mathrm{clsd}}^{2n})). \label{2.5}
\end{eqnarray}
We call them the closed characters. The restriction map to
Spec$(k(X))$ fulfills

\begin{eqnarray}
  AD^1(X)_{\mathrm{clsd}} & \stackrel{\cong}{\rightarrow} & H^0(X, \Omega_{X/k, \mathrm{clsd}}^1/d \mathrm{\log}\mathcal
  {O}_X^{\times})\subset AD^1(k(X))_{\mathrm{clsd}} \nonumber \\
  AD^n(X)_{\mathrm{clsd}} & \rightarrow & H^0(X, \mathcal {H}_{DR}^{2n-1})\subset H_{DR}^{2n-1}(k(X)/k), \quad
  \mbox{for } n\geq2 \label{2.6}
\end{eqnarray}
(see [2]). Here $\mathcal {H}_{DR}^p$ is the Zariski sheaf of de
Rham cohomology.

If $k$ is the field of complex numbers $\mathbb{C}$, one can
change from the Zariski topology to the analytic one. This yields
a map
\begin{eqnarray}
  AD(X) & = & \oplus_nAD^n(X)\stackrel{l}{\rightarrow} \nonumber \\
  D(X) & = & \oplus_nD^n(X)=\oplus_n\mathbb{H}^{2n}(X_{\mathrm{an}}, \mathbb{Z}(n)\rightarrow
   \mathcal {A}_{X}^0\stackrel{d}
  {\rightarrow}\ldots \rightarrow \mathcal {A}_X^{2n-1}). \label{2.7}
\end{eqnarray}
Here $D(X)$ is the ring of differential characters defined by
Cheeger-Simons ([10]). One has
\begin{equation} \label{2.8}
  \iota(AD(X)_{\mathrm{clsd}})\subset\oplus_nH^{2n-1}(X_{\mathrm{an}}, \mathbb{C}/\mathbb{Z}(n))\subset D(X).
\end{equation}

It is classical that $AD^1(X)$ is the group of isomorphism classes
of lines bundles $L$ with connection $\nabla$. If
${\xi}_{ij}\in\Gamma(U_{ij}, \mathcal {O}_X^\times)$ is a cocycle
of $L$ in a local frame $e_i$ of $L$ on $U_i$, and
$\nabla(e_i)=\alpha_i\in\Gamma(U_i, \Omega_{X/k}^{1})$ is the
local form of the connection, then
$d\log\xi_{ij}=\alpha_j-\alpha_i=\delta(\alpha)_{ij}$ defines the
Cech cocycle of $c_1((L, \nabla))$. In [12], [13], [15], we
generalize this class.\vskip 2mm

\noindent {\bf Theorem 2.1} ([12], [13], [15]){\bf .} \it
Associated to an algebraic bundle E with  connection $\nabla$
{\rm(}resp. with integrable connection{\rm)}, one has
characteristic classes $c_n((E, \nabla)) \in AD^n (X)$ {\rm
(}resp. $\in AD^n(X)_{\rm clsd})$. These classes satisfy the
following properties:
\begin{enumerate}
\item[{\rm (1)}] The classes $c_{\ast}((E,
\mathrm{\nabla}))\in AD^\ast(X)$ are functorial and additive.

\item[{\rm (2)}] $c_1((E, \mathrm{\nabla}))$ is the isomorphism class of
$(\det(E), \det({\mathrm{\nabla}}))$.

\item[{\rm (3)}] forget $(c_n((E, \mathrm{\nabla})))=c_n(E)\in CH^n(X)$
is the Chern class of the underlying algebraic bundle $E$ in the
Chow group.

\item[{\rm (4)}] $d(c_n((E, \mathrm{\nabla})))=c_n(\mathrm{\nabla}^2)\in
H^0(X, \Omega_{X/k, \mathrm{clsd}}^{2n})$ is the Chern-Weil form
which is the evaluation of the invariant polynomial $c_n$ on the
curvature $\mathrm{\nabla}^2$.

\item[{\rm (5)}] The restriction to the generic point $c_n((E,
\mathrm{\nabla})|_{k(X)})$ is the algebraic Chern-Simons invariant
$CS_n ((E, \mathrm{\nabla}))$ defined in {\rm[4]}.\  It has values
in $H^0 (X, $  \linebreak  $\Omega_{X/k}^1/d\log\mathcal
{O}^{\times})$\ $(resp. H^0(X, \Omega_{X/k, {\rm
clsd}}^1/d\log\mathcal {O}^{\times}))$ for $n=1$, and in $H^0 (X,
\linebreak \Omega_{X/k}^{2n-1}/d\Omega_{X/k}^{2n-2})\ (resp. H^0
(X, \mathcal{H}_{DR}^{2n-1}))$ for $n\geq2$.

\item[{\rm (6)}] If $k\subset\mathbb{C}$, then $\iota (C_n ((E,
\mathrm{\nabla})))\in D^n(X) (resp. \in H^{2n-1}(X_{\mathrm{an}},
\mathbb{C}/\mathbb{Z}(n)))$ is the differential character defined
by Chern-Cheeger-Simons, denoted by
$c_n(E_{\mathrm{an}}^{\mathrm{\nabla}})\in
H^{2n-1}(X_{\mathrm{an}}, \mathbb{C}/\mathbb{Z}(n))$ if
$\mathrm{\nabla}$ is integrable.
\end{enumerate} \rm

If $k=\mathbb{C}, \mathrm{\nabla}$ is flat and the underlying
monodromy is finite, then the existence of $c_n{((E,
\mathrm{\nabla}))} $immediately implies that the Chern-Simons
classes in $H^{2n-1}(X_{\mathrm{an}}, \linebreak
\mathbb{Q}(n)/\mathbb{Z}(n))$ are in the smallest possible level
of the coniveau filtration ([14]).

If $X$ is complex projective smooth, $\mathrm{\nabla}$ is
integrable and $n\geq2$,  we relate $CS_n((E, \mathrm{\nabla}))$
for $n\geq2$ to the (generalized) Griffiths'\, group
$\mathrm{Griff}^n(X)$. It consists of cycles which are homologous
to $0$ modulo those which are homologous to $0$ on some
divisor([4], definition 5.1.1). For $n=2, [2]$ implies that
$\mathrm{Griff}^2(X)$ is the classical Griffiths'\, group. For
$n\geq2$, Reznikov's theorem  ([19]) (answering positively Bloch's
conjecture [3]), together with the existence of the lifting $c_n
((E, {\mathrm{\nabla}}))$, imply that the classes $CS_n((E,
\mathrm{\nabla}))$ lie in the image Im of the global cohomology
$H^{2n-1}(X, \mathbb{Q}(n))$ in $H^0(X, \mathcal
{H}^{2n}(\mathbb{Q}(n))$. This subgroup Im maps to
$\mathrm{Griff}^{n}$(X). One has

\noindent {\bf Theorem 2.2} ([4], Theorem 5.6.2){\bf .} \it
$$ \mbox{\rm image of\  }\ CS_n((E, \nabla))\in {\rm{Griff}}^n(X)\otimes \mathbb{Q}$$
is the Chern class $c_n^{\rm{Griff}}(E)\otimes \mathbb{Q}$ of the
underlying algebraic bundle E. Moreover, $CS_n((E,
\mathrm{\nabla}))=0$ if and only if \
$c_n^{\mathrm{Griff}}(E)\otimes \mathbb{Q}=0$. \rm

A relative version $AD(X/S)$ of $AD(X)$ is defined in [6]. We give here an example of application.

\noindent {\bf {Theorem 2.3}}([6], Corollary 3.15){\bf .} \it Let
$f:X\rightarrow S$ be a smooth projective family of curves over a
field $k$. Let $(E, \mathrm{\nabla}_{X/S})$ be a bundle with a
relative connection. Then there are classes $c_2((E,
\mathrm{\nabla}))\in AD^2(X/S):=\mathbb{H}^2(X, \mathcal{K}_2
\stackrel{d\log} \longrightarrow \Omega_S^1\otimes_{X/S}^1 )$\
lifting the classes $c_2(E)\in CH^2(X)$. There is a trace map
$f_\ast: AD^2(X/S)\rightarrow AD^1(S)$ compatible with the trace
map on Chow groups $f_\ast:CH^2(X)\rightarrow CH^1(S)$. Thus
$f_\ast c_2(E, \mathrm{\nabla})$ is a connection on the line
bundle $f_\ast c_2(E)$, which depends functorially on the choice
of $\mathrm{\nabla}_{X/S}$ on $E$. \rm

We now study the behavior at $\infty$ of $CS_n((E,
\mathrm{\nabla}))$ for $n\geq2$ in characteristic $0$. Let
$j:X\rightarrow \bar{X}$ be a smooth compactification of $X$.
Recall that a de Rham class $\in H_{DR}^q(k(X)/k)$ at the generic
point is called $unramified$ if it lies in $H^0(\bar{X}, \mathcal
{H}_{DR}^q\subset H_{DR}^q(k(X)/k) ([2])$.

\noindent {\bf Theorem 2.4.} \it We assume k to be of
characteristic $0$, and $\mathrm{\nabla}$ to be integrable. Then
$CS_n((E, \mathrm{\nabla}))$ is unramified for $n\geq 2$. \rm

\noindent {\bf Proof.}\ If $(E, \mathrm{\nabla})$ is regular
singular, this is shown in [4], theorem 6.1.1. In general, one may
argue as follows. One has
$$H^0(\bar{X}, \mathcal {H}_{DR}^{2n-1})=\mathrm{Ker}(H^0(X,
\mathcal {H}_{DR}^{2n-1})\stackrel{\mathrm{Res}}{\longrightarrow}
\oplus_{x\in X^{(1)}}H_{DR}^{2n-2}(x/k)).$$ Thus it is enough to
show that the residue map at each generic point at $\infty$ dies.
At a smooth point of a divisor $D$ at $\infty$, the residue
depends only on the formal completion of $X$ along $D$. So we may
assume that $\mathrm{\nabla}$ is a connection on $\mathcal
{O}_X=k(D)[[x]]$, integrable relative to $k$. By a variant (see
[1], proposition 5.10.) of Levelt's theorem ([17]) for absolute
flat connections, there are a finite extension $K\supset k(D)$,
and a ramified extension $\pi:K[[x]]\subset K[[y]], y^N=x$ for
some $N\in \mathbb{N} \setminus \{0\}$ such that
$\pi^\ast(\mathrm{\nabla})=\oplus(L\otimes U)$. Here $L$ is
integrable of rank 1 and $U$ is integrable with logarithmic poles
along $y$=0. Since
$\mathrm{Res}_{y=0}\pi^{\ast}(\alpha)=N\pi^{\ast}(\mathrm{Res}_{x=0}(\alpha))$,
and $H_{DR}^p(k(D)/k)\subset H_{DR}^p(K/k)$, functoriality and
additivity of the classes imply that we may assume
$\mathrm{\nabla}=L\otimes U$ on $K((x))$ with $K=k(D)$, In a local
frame we have the equations  $U=\Gamma\dfrac{dx}{x}+\Sigma$ where
$\Gamma\in GL(r, K[[x]]), \Sigma\in M(r, K[[x]])\otimes
\Omega_K^1$, and $L=d(f)+\lambda\dfrac{dx}{x}+\beta$ where $f\in
K((x)), \lambda\in k, \beta\in \Omega_K^1$ with $d\beta=0$. The
explicit formula Res Tr($A(d(A))^{n-1})\in H_{DR}^{2n-1}(K)$ of
$CS_n((E, \mathrm{\nabla}))\in H_{DR}^{2n-1}(K((x)))$ and [4],
prop. 5.10 in the logarithmic case, imply that
$$CS_n((E, \mathrm{\nabla}))=\mathrm{Tr}(d(f)+\lambda\dfrac{dx}{x}+\beta)(d(\Gamma)\dfrac{dx}{x}+
d\Sigma)^{n-1}.$$ This is the sum of 2 terms with rational
coefficients, Res Tr$(d(f)+\beta)(d(\Sigma)^{n-2}$
$d(\Gamma)\dfrac{dx}{x})$ and Res
Tr$(\lambda\dfrac{dx}{x}(d(\Sigma)^{n-1}))$. Both terms are
obviously exact.

\noindent {\bf {Discussion 2.5.} \rm  We assume here
$k=\mathbb{C}$, $\mathrm{\nabla}$ is integrable and $n\geq 2$. We
consider the image
$c_n(E_{\mathrm{an}}^{\mathrm{\nabla}}|_{\mathbb{C}(X)})\in
H^0(\bar{X}, {\mathcal H}^{2n-1}(\mathbb{C}/\mathbb{Z}(n)))$ of
$CS_n((E, \mathrm{\nabla}))\in H^0(\bar{X}, \mathcal
{H}_{DR}^{2n-1})$. When $X$ is not compact, there is Deligne's
unique algebraic ($E, \mathrm{\nabla}$) with regular singularities
at $\infty$ with the given underlying local system
$E_{an}^{\mathrm{\nabla}}$([11]), but there are many irregular
connections ($E, \mathrm{\nabla}$). The topological class
$C_n(E_{\mathrm{an}}^{\mathrm{\nabla}})\in
H^{2n-1}(X_{\mathrm{an}}, \mathbb{C}/\mathbb{Z}(n))$ is not, a
priori, extendable to $\bar{X}$, but we have seen that its
restriction to Spec$(\mathbb{C}(X))$ is unramified. \rm

There is on $X$ a fundamental system of Artin neighborhoods $U$
which are geometrically successive fiberings in affine curves.
Topologically they are $K\pi_1$ and their fundamental group is a
successive extension of free groups in finitely many letters. On
such an open $U$, the class $c_n(E_{\mathrm{an}}^\mathrm{\nabla})$
lies in $H^{2n-1}(U_{\mathrm{an}},
\mathbb{C}/\mathbb{Z}(n))=H^{2n-1} (\pi_1(U_{\mathrm{an}}, u),
\mathbb{C}/\mathbb{Z}(n))$.

If $U$ is such that $\pi_1(U_{\mathrm{an}}, u)$ is isomorphic as
an abstract group to $\pi_1(V_{\mathrm{an}}, v)$, where $V$ is an
Artin neighborhood on a rational variety, then
$E_{\mathrm{an}}^{\mathrm{\nabla}}|_U$ becomes a representation of
$\pi_1(V_{\mathrm{an}}, v)$, and since then $H^0(\bar{V}, \mathcal
{H}_{DR}^{2n-1})=0$ for $n\geq 2$ and $V\subset \bar{V}$ a good
compactification, one obtains
$c_n(E_{\mathrm{an}}^{\mathrm{\nabla}}|_{\mathbb{C}(X)})=0, n\geq
2$ in this case. Such an example is provided by a product of
smooth affine curves of any genus. It has the fundamental group of
a product of $\mathbb{P}^1$ minus finitely many points.

\noindent {\bf Question 2.6.} In view of the previous discussion,
we may ask what complex smooth varieties $X$ are dominated by
$h:Y\rightarrow X$ proper, with $Y$ smooth, such that $Y$ has an
Artin neighborhood, the fundamental group of which is the
fundamental group of an Artin neighborhood on a rational variety,
or more generally of a variety for which $H^0$(smooth
compactification, ${\mathcal {H}}_{DR}^n)=0$ for $n\geq 1$. We
have seen that this would imply vanishing modulo torsion of
$c_n(E_{\mathrm{an}}^{\mathrm{\nabla}}|_{\mathbb{C}(X)})$, $n\geq
2$, or equivalently $CS_n((E, \mathrm{\nabla}))\in H^0(\bar{X},
{\mathcal H}^{2n-1}(\mathbb{Q}(n)))$.

On the other hand, if $X$ is projective smooth, Reznikov's theorem
([19]) shows vanishing modulo torsion of the Chern-Cheeger-Simons
classes in $H^{2n-1}(X_{\mathrm{an}}$, $\mathbb{C}/
\mathbb{Q}(n))$. It is a consequence of Simpson's nonabelian Hodge
theory on smooth projective varieties. Our classes $CS_n((E,
\mathrm{\nabla}))$ live at the generic point of $X$. We don't have
a nonabelian mixed Hodge theory at disposal. Yet one may ask
whether it is always true that $CS_n((E, \mathrm{\nabla}))\in
H^0(\bar{X}, {\mathcal H}^{2n-1}(\mathbb{Q}(n)))$  for $n\geq 2$,
even if many $X$ don't have the topological property explained
above.

\section{The behavior of\  the \ algebraic \ Chern-Simons \\
 classes in families in the regular singular case}
\label{section 3} \setzero

\vskip-5mm \hspace{5mm}

The algebraic Chern-Simons invariants $CS_n((E, \mathrm{\nabla}))$
have been studied in a family in [5]. Given $f:X\rightarrow S$ a
proper smooth family, and ($E, \mathrm{\nabla}$) a flat connection
on $X$, the Gau\ss-Manin bundles
$$R^if_\ast(\Omega_{X/S}^{\bullet}\otimes E, \mathrm{\nabla}_{X/S})$$
carry the Gau\ss-Manin connection $GM^i(\mathrm{\nabla})$. We give
a formula for the invariants $CS_n((GM^i(\mathrm{\nabla})-$rank
$(\nabla)\cdot GM^i(d)))$ on $S$, as a function of $CS_n((E,
\mathrm{\nabla}))$ and of characteristic classes of $f$. Here ($O,
d$) is the trivial connection.

More generally, we may assume that $f$ is smooth away from a
normal crossings divisor $T\subset S$ such that
$Y=f^{-1}(T)\subset X$ is a normal crossings divisor with the
property that $\Omega_{X/S}^1(\log Y)$ is locally free. Then ($E,
\mathrm{\nabla}$) has logarithmic poles along $Y\cup Z$ where
$Y\cup Z\subset X$ is a normal crossings divisor, still with the
property that $\Omega_{X/S}^1(\log(Y+Z))$ is locally free. That is
$Z$ is the horizontal divisor of singularities of
$\mathrm{\nabla}$. The formula involves the top Chern class
$c_d(\Omega_{X/S}^1(\log (Y+Z))\in \mathbb{H}^d(X, \mathcal
{K}_d\rightarrow\oplus_i\mathcal {K}_{Z_{i, d}})$, rigidified by
the residue maps $\Omega_{X/S}^1(\log (Y+Z))\rightarrow \mathcal
{O}_{Z_i}$, as defined by T. Saito in [20]. One of its main
features is that $CS_n((GM^i(\mathrm{\nabla})-$rank
$(\mathrm{\nabla})\cdot GM^i(d)))$ vanishes if $CS_n((E,
\mathrm{\nabla}))$ vanishes. It is

\noindent {\bf Theorem 3.1} ([5], Theorem 0.1){\bf .}
\begin{eqnarray*}
  & & CS_n(\sum(-1)^i(GM^i(\mathrm{\nabla})-{\rm rank}(\mathrm{\nabla})\cdot GM^i(d))) \\
  &=&(-1){^{\mathrm{dim}(X/S)}}f_{\ast}{c_{\mathrm{dim}(X/S)}}(\Omega_{X/S}^1(\log
  (Y+Z)), {\rm res})\cdot CS_n((E, \mathrm{\nabla})).
\end{eqnarray*}
Here $\cdot$ is the cup product of the algebraic Chern-Simons invariants with this rigidified class, which is well
defined, as well as the trace $f_\ast$ to $S$.

\noindent {\bf Discussion 3.2.}\ One weak point of the method used
in [5] is that it does not allow to understand a formula for the
whole invariants $c_n((E, \mathrm{\nabla}))$, but only for
$CS_n(E, \mathrm{\nabla})$. Indeed, we use the explicit formula
studied in [4] to compute it, which can't exist for the whole
class in $AD(X)$, as it in particular involves the Chern classes
of the underlying algebraic bundle $E$ in the Chow group. \rm

\section{ The determinant of the Gauss-Manin connection: the irregular rank 1 case} \label{section 4}\setzero

\vskip-5mm \hspace{5mm}

Now we no longer assume that ($E, \mathrm{\nabla}$) is regular
singular at $\infty$. In the next two sections, we reduce
ourselves to the case where $f:X\rightarrow S$ is a family of
curves, and we consider only the determinant of the Gau\ss-Manin
connection. That is we consider
\begin{eqnarray*}
  \det(GM)&:=&\sum_i(-1)^ic_1(GM^i) \\
  & \in & \mathbb{H}^1(S, {\cal{O}}_S^{\times} \stackrel{d\mathrm\log}{\longrightarrow}\Omega_S^1)\subset
  \Omega_{k(S)}^1/d\log(k(S))^{\times}.
\end{eqnarray*}
Since the determinant is recognized at the generic point of $S$,
we replace $S$ by its function field $K:=k(S)$ in the next two
sections. In other words, $X/K$ is an affine curve. Let
$\bar{X}/K$ be its smooth compactification.

In this section, we assume that the integrable connection ($L,
\mathrm{\nabla}$) we start with on $X$ has rank 1. Following
Deligne's idea (see his 1974 letter to Serre published as an
appendix to [7], we first reduce the problem of computing the
determinant of the cohomology of $\mathrm{\nabla}$ on the curve to
the one of computing the determinant of cohomology of an
integrable invariant connection, still denoted by
$\mathrm{\nabla}$, on a generalized Jacobian. More precisely,
$\nabla$ has a divisor (with multiplicities) of irregularity
$\sum_im_ip_i$, where $m_i-1$=irregularity of $\nabla$ in $p_i\in
\bar{X}\setminus X$. On the Jacobian $G=\mathrm{Pic}(\bar{X},
\sum_im_ip_i)$ of line bundles trivialized at the order $m_i$ at
the points $p_i$, there is an invariant connection which pulls
back to $\mathrm{\nabla}$ via the cycle map. On the torsor
$p^{-1}(\omega_{\bar{X}}(\sum_im_ip_i))$ under the affine group
$p^{-1}({\cal{O}}_{\bar{X}})$, where $p:G\rightarrow
\mathrm{Pic}(\bar{X})$, one considers the hypersurface
$\Sigma:\sum_i \mathrm{res}_{p_i}=0$. We show that the relative
invariant connection $\mathrm{\nabla}_{/K}$ on $G$, while
restricted to $\Sigma \subset
p^{-1}(\omega_{\bar{X}}(\sum_im_ip_i))$, acquires exactly one zero
which is a $K$-rational point $\kappa$ of $G$. Restricting
$\mathrm{\nabla}$ to this special point yields a connection
$\mathrm{\nabla}|_{\kappa}$ on $K$. The formula then says that the
determinant of the Gau\ss-Manin connection is the sum of this
connection $\mathrm{\nabla}|_{\kappa}$ and of a 2-torsion term,
which we describe now. In a given frame of $L$ at a singularity
$p_i$, the local equation of the connection is
$\alpha_i=a_i\dfrac{dt_i}{t_{i}^{m_i}}+$lower order terms, where
$t_i$ is a local parameter and $a_i\in k^{\times}$. Then the
2-torsion connection $\frac{m_i}{2}d\log a_i\in
\Omega_{K/k}^1/d\log K^{\times}$ does not depend on the choices
and is well defined. The 2-torsion term is the sum over the
irregular points of these 2-torsion connections. Summarizing, one
has

\noindent {\bf Theorem 4.1} ([7], Theorem 1.1){\bf .}
\begin{eqnarray*}
  & & \det \left( \sum_i(-1)^iH^i(X, (\Omega_{X/K}^{\bullet} \otimes L, \mathrm{\nabla}_{X/K}), GM^i(\mathrm{\nabla}) \right) \\
  & = & (-1)\mathrm{\nabla}|_{\kappa}+\sum_{i, m_i\geq2}\dfrac{m_i}{2}d \log a_i\in\Omega_{K/k}^1/d \log K^\times.
\end{eqnarray*}

\noindent {\bf Discussion 4.2.}\ The formula described above is
global. As such, it has a spirit which is different from Deligne's
formula describing the global $\epsilon$-factor of an $\ell$-adic
character over a curve over a finite field as a product of local
$\epsilon$-factors. However, choosing another $K$-rational point
$\kappa'\in p^{-1}(\omega_{\bar{X}}(\sum_im_ip_i))$, it is easy to
write the difference
$\mathrm{\nabla}|_{\kappa'}-\mathrm{\nabla}|_{\kappa}$ as a sum of
explicitly given connections on $K$. One has
$u\cdot\kappa'=\kappa$, where $u=\prod u_i\in
p^{-1}\mathcal{O}_{\bar{X}}=\prod_i ({\cal{O}}_{\bar{X},
p_i}/{\frak{m}}_{\frak{p_i}})^{\times}/K^{\times}$ . Then
$\mathrm{\nabla}_{\kappa'}=\mathrm{\nabla}|_{\kappa}-\sum_i\mathrm{Res}_{p_i}d\log
u_i\wedge\alpha_i$. Correspondingly, one may write the right hand
side of the formula above as
$$(-1)\mathrm{\nabla}|_{\kappa'}+\sum_i\biggl(\sup\big(1, \dfrac{m_i}{2}\big)d\log
a_i+\mathrm{Res}_{p_i}d\log u_i\wedge\alpha_i\biggr). $$ In
particular, the choice of some differential form
$\nu\in\omega_{\bar{X}/K}\otimes K(X)$, generating
$\omega_{\bar{X}/K}(\sum m_ip_i)$ at the point $p_i$, defines a
trivialization of $\omega_{\bar{X}/K}(\sum_im_ip_i)$ thus a point
$\kappa(\nu)$. We write $\alpha_i=g_i\nu$ with $g_i\in({\cal
{O}}_{X, p_i}/{\frak m}_{\frak p_i})^{\times}$. Then  the formula
reads

\noindent {\bf Theorem 4.3} ([7], Formula 5.4){\bf .}
\begin{eqnarray*}
  & & \det \left(\sum_i(-1)^i(H^i(X, (\Omega_{X/K}^{\bullet}\otimes L, \mathrm{\nabla}_{X/K}), GM^i(\mathrm{\nabla})
  \right) \\
  & = & (-1)\det (\mathrm{\nabla})|_{\kappa(\nu)}+\sum_i\left(\sup\big(1, \dfrac{m_i}{2}\big)d\log
  a_i+\mathrm{Res}_{p_i}dg_ig_i^{-1}\wedge\alpha_i\right).
\end{eqnarray*}

\section{The determinant of the Gauss-Manin connection: the irregular higher rank case} \label{section 1}\setzero

\vskip-5mm \hspace{5mm}

We assume in this section that we have an affine curve $X$ over
$K=k(S)$, $k$ of characteristic zero as in the rank 1 case. The
integrable connection $(E, \mathrm{\nabla})$ we are given on $X$
has higher rank $r$.

In the rank one case, for any rank one bundle contained in
$j_{\ast}E$, the equation of the connection in a local formal
frame at a singular point is of the shape
$\alpha=a\dfrac{dt}{t^m}+\dfrac{\beta}{t^{m-1}}$, where
$m\in\mathbb{N}, a\in K[[t]]^{\times}$, $\beta\in\Omega_{K}^1
\otimes K[[t]]$. In particular, $(m-1)$ is the irregularity of the
connection ([11]). Here $t$ is a local parameter. If $r>1$, it is
no longer the case that $j_{\ast}E$ necessarily contains a rank
$r$ bundle such that in a local formal frame of this bundle, the
local equation has the shape
$A_i=a_i\dfrac{dt_i}{t_i^m}+\dfrac{\beta_i}{t_i^{m-1}}$, with
$a_i\in GL(r, K[[t_i]]), \beta_i\in\Omega_K^1 \otimes M(r,
K[[t_i]])$. We call an integrable connection $(E,
\mathrm{\nabla})$ with this existence property an $admissible$
connection.

Even if $(E, \nabla)$ is admissible, its determinant connection
$\det(E, \mathrm{\nabla})$ might have much lower order poles (for
example trivial). This indicates that one can not extend directly
in this form the formula 4.1. However, assuming $(E,
\mathrm{\nabla})$ to be admissible and choosing some
$\nu\in\omega_{\bar{X}/K}\otimes K(X)$ which generates
$\omega(\sum_im_ip_i)$ at $p_i$ as for formula 4.3, the right hand
side of 4.3 makes sense, if one replaces $d \log a_i$ by
\linebreak $d \log\det(a_i)$. Using global methods inspired by the
Higgs correspondence between Higgs fields and connections on
complex smooth projective varieties ([21]), one is able to prove
the ``same" formula as 4.3 in the higher rank case on
$\mathbb{P}^1$.

\noindent {\bf Theorem 5.1} ([8], Theorem 1.3){\bf .} \it If $(E,
\mathrm{\nabla})$ is admissible and has at least one irregular
point, and if $\nu\in\omega_{\bar{X}/K}\otimes K(X)$ generates
$\omega(\sum_im_ip_i)$ at the points $p_i$, then
\begin{eqnarray*}
  & & \det \left( \sum_i(-1)^i(H^i(X, (\Omega_{X/K}^{\bullet}\otimes L,  \mathrm{\nabla}_{X/K}), GM^i(\mathrm{\nabla}))
  \right) \\
  & = & (-1)\mathrm{\nabla}|_{\kappa(\nu)}+\sum_i\left(\sup\big(1, \dfrac{m_i}{2}\big)d\log
  \det(a_i(p_i))+\mathrm{Tr \ Res}_{p_i}dg_ig_i^{-1}\wedge A_i\right).
\end{eqnarray*}
\rm

The connection Res $\mathrm{Tr}_{p_i}dg_ig_i^{-1}\wedge
A_i\in\Omega_K^1/d\log K^{\times}$ is well defined, as well as the
2-torsion connection $\sup\big(1, \dfrac{m_i}{2}\big)d\log\det
(a_i(p_i))$.

However, one needs a different method in order to understand the
contribution of singularities in which $(E, \mathrm{\nabla})$ is
not admissible.

We describe now the origin of the method contained [1]. It is
based on the idea that Tate's method ([22]) applies for
connections.

Locally formally over the Laurent series field $K((t)), E$ becomes
a $r$-dimensional vector space over $K((t))$. The relative
connection
$\mathrm{\nabla}_{K((T))/K}:E\rightarrow\omega_{K((t))/K}\otimes
E$ is a Fredholm operator. This means that
$H^i(\mathrm{\nabla}_{X/K}), i=0, 1$ are finite dimensional
$K$-vector spaces, and that $\mathrm{\nabla}_{X/K}$ carries
compact lattices to compact lattices. Let $E\cong\oplus_1^rK((t))$
be the choice of a local frame. A compact lattice is a
$K$-subspace of $E$ which is commensurable to $\oplus_1^rK[[t]]$.
Given $0\neq \nu \rightarrow \in \omega_{K((t))/K}$, one composes
$\mathrm{\nabla}_{K((t))/K,
\nu}:=\nu^{-1}\circ\mathrm{\nabla}_{K((t))/K}:E\rightarrow E$ to
obtain a Fredholm endomorphism. To a Fredholm endomorphism
$A:E\rightarrow E$, one associates a \linebreak 1-dimensional
$K$-vector space $\lambda(A)=\det(H^0(A))\otimes\det(H^1(A))^{-1}$
together with the degree $\chi(A)=\dim H^0(A)-\mathrm{dim}H^1(A)$.
We call this a {\em super-line}. It does not refer to the topology
defined by compact lattices. Then one measures how $A$ moves a
compact lattice $L\subset E$. First for 2 lattices $L$ and $L'$,
one takes a smaller compact lattice $N\subset L \cap L'$ and
defines $\det(L:L'):=\det(L/N)\cdot\det(L'/N)^{-1}$, where $\cdot$
is the tensor product of super-lines and $\det(L/N)$ has degree
dim($L/N$). This does not depend on the choice of $N$. Then one
defines asymptotic superlines. The compact one is
$\lambda_c(A)=\det(A(L):L)\cdot\det(L\cap \mathrm{Ker}(A))$ and
the discrete one is
$\lambda_d(A)=\det(L:A^{-1}(L))\cdot\det(V/(L+A(V))$. They do not
depend on the choice of $L$. Taking
$0\neq\nu\in\omega_{\bar{X}/K}\otimes K(X)$ a rational
differential form, and $E_{\mathrm{min}}$ the minimal extension of
$E, X$ the complement of the singularities of $\nu$ and
$\mathrm{\nabla}_{X/K}$, one has
$$\det H^{\ast}(X/K, E_{\mathrm{min}})=\otimes_{X\in\bar{X}\setminus
X}\lambda_d(\mathrm{\nabla}_{K((t))/K,
\nu})\cdot\det(E^{\mathrm{\nabla}})^{-1}$$ as a product of
discrete lines.

On the other hand, one has the relation
$\lambda(A)=\lambda_d(A)\cdot\lambda_c(A)^{-1}$. One easily
computes that $\chi(\mathrm{\nabla}_{K((t))/K, \nu})=0,
\lambda(\mathrm{\nabla}_{K((t))/K, \nu})=1$. Setting
$$\epsilon_x(\mathrm{\nabla}_{K((t))/K, \nu})=\lambda_c(\mathrm{\nabla}_{K((t)/K, \nu})
\cdot\det(E^{\mathrm{\nabla}})^{-1}, $$ this implies immediately
the product formula.

\noindent {\bf Theorem 5.2} ([1], (1.3.1)){\bf .}
$$ \det H^{\ast}(X/K, E_{\mathrm{min}})=\oplus_{x\in\bar{X}
\setminus X}\epsilon_x(\mathrm{\nabla}_{K((t))/K, \nu}).$$

It remains to endow the local $\epsilon$ lines with a connection,
compatible with the Gau\ss-Manin connection on the left. One
chooses a section of the vector fields $T_{K/k} \subset T_{X/k}$
and a relative differential form $\nu$ which is annihilated by
this section. One applies Grothendieck's definition of a
connection. The Gau\ss-Manin connection is given by the
infinitesimal automorphism $p_1^{\ast}\mathrm{det}H^*(X/K,
E_{\mathrm{min}})\rightarrow p_2^{\ast}\det H^{\ast}(X/K,$ $
E_{\mathrm{min}})$ on $K\otimes_k K$, induced by $\tau:
p_1^{\ast}E_{\mathrm{min}}\rightarrow p_2^{\ast}E_{\mathrm{min}},
\tau \in T_{K/k}$, which by the choice commutes to
$\mathrm{\nabla}_{X/K}$. By functoriality of the $\epsilon$ lines,
this defines a $K\otimes_kK$ homomorphism
$p_1^{\ast}\epsilon\rightarrow p_2^{\ast}\epsilon$. This is the
$\epsilon$-connection.

\noindent {\bf Theorem 5.3} ([1], Theorem 5.6){\bf .} \it For an
admissible connection of local equation
$A=a\dfrac{dt}{t^m}+\dfrac{\beta}{t^{m-1}}$, with $m\geq2$, the
local $\epsilon$-connection is
\begin{eqnarray*}
  \epsilon\big(\mathrm{\nabla}_{K((t))/K}, \dfrac{dt}{t^m}\big) & = & \mathrm{Tr} \mathrm{Res}_{t=0}daa^{-1}A+
  \dfrac{m}{2}d\log\det(a(t=0)) \\
  & \in & \Omega_{K/k}^1/d\log K^{\times}.
\end{eqnarray*}
\rm

The restriction on the choice of $\nu$ given by the commutativity constraint with some lifting of vector fields of
$K$ is not necessary. The construction is more general.

Given a relative connection $\mathrm{\nabla}_{K((t))/K}$, the
$\epsilon$-lines for $0\neq\nu\in\omega_{K((t))/K}^{\times}$ build
a super-line bundle on the ind-scheme
$\omega_{K((t))/K}^{\times}$. The line bundle obeys a connection
relative to $K$ on $\omega_{K((t))/K}^{\times}$. Formula 5.2
identifies line bundles with connections relative to $K$, where
the left hand side carries the constant connection. The choice of
an integrable lifting $\mathrm{\nabla}$ of $\mathrm{\nabla}_{X/K}$
yields a lifting of the relative connection on the $\epsilon$ line
to an integrable connection relative to $k$. Formula 5.2
identifies line bundles with integrable connections where the left
hand side carries the Gau\ss-Manin connection ([1], 1.3).

The $\epsilon$ lines and connections are additive in exact sequences and compatible with push-downs. By a variant
of Levelt's theorem for integrable formal connections, this allows to show that all connections are induced from
admissible ones, for which we have the formula 5.3.

\noindent {\bf Question 5.4.}\  We don't know how to precisely
relate the algebraic group viewpoint developed to treat the rank 1
case, and the special rational point found there, with the
polarized Fredholm line method which works in general.

\noindent {\bf Acknowledgements.}\  I thank the mathematicians I
have worked with on the material exposed in those notes. A large
part of it has been jointly developed with Spencer Bloch. It is a
pleasure to acknowledge the impact of his ideas on a programme I
had started earlier and we continued together. I thank Alexander
Beilinson. An unpublished manuscript of his and David Kazhdan
allowed me to understand completely one of the two constructions
explained in [15]. His deep viewpoint reflected in [1] changed the
understanding of the formula we had as explained in [8]. I thank
Pierre Deligne, whose ideas on epsilon factors have shaped much of
my thinking. His letter to Serre on the rank 1 case is published
as an appendix to [7], but the content of his seminar at the IHES
in 1984 has not been available to me. I thank Takeshi Saito for
his willingness to explain different aspects of the $\ell$ -adic
theory.

\label{lastpage}

\end{document}